\documentclass[letterpaper]{amsart}
\usepackage{fullpage}

\usepackage{amssymb,amsfonts,amscd,verbatim}
\usepackage[leqno]{amsmath} 
\usepackage{mathrsfs}
\usepackage{xy} % for diagrams
\usepackage[normalem]{ulem}% for double underline but \emph in italics
\makeatletter

\newcommand{\id}{\mathrm{id}}

\newcommand{\fl}{{\rm  fl}}
\newcommand{\ft}{{\rm  ft}}
\newcommand{\tor}{{\rm  tor}}
\newcommand{\fr}{{\rm  fr}}
\newcommand{\Z}{{\mathbb Z}}
\newcommand{\Q}{{\mathbb Q}}
\newcommand{\G}{{\mathbb G}}
\newcommand{\Gr}{{\rm Gr}}
\newcommand{\GGr}{{\bf Gr}}
\newcommand{\spec}{{\rm Spec}}
\newcommand{\Pic}{{\rm Pic}}

\newtheorem{lemma}{Lemma}{\bf}{\it}
\newtheorem{theorem}[lemma]{Theorem}{\bf}{\it}
\newtheorem{proposition}[lemma]{Proposition}{\bf}{\it}
{\bf}{\rm}
\newtheorem{remark}[lemma]{Remark}{\bf}{\rm}
{\bf}{\it}
{\bf}{\it}

\begin{document}
\baselineskip 15pt

\input xy     %%for diagrams
\xyoption{all}%%

\title{On torsors under abelian varieties}
\author{A. Bertapelle}
\date{July 26, 2011}

\begin{abstract} Let $A_K$ be an abelian variety over a local field $K$ of mixed characteristic and with algebraically closed residue field. We provide a geometric construction (via the relative Picard functor) of the Shafarevich duality  between the group of  isomorphism classes of torsors under $A_K$ and the ``fundamental group'' of the N\'eron model of the dual abelian variety $A_K'$. An analogous construction works over  fields of positive characteristic $p$ providing a duality on the prime-to-$p$ parts.
\end{abstract}
\maketitle

Let $R$ be a  complete discrete valuation ring
with field of fractions $K$ and algebraically closed residue field $k$.
Let us denote by $j\colon \spec(K)\to \spec(R)$ the usual open immersion.
Let $A_K,A^\prime_K$ be dual abelian varieties over $K$, and $A,A^\prime$ their N\'eron models. 

Shafarevich's pairing (cf.  \cite{Beg}, \cite{Bes}, \cite{Ber}) provides an isomorphism
\begin{equation}\label{eq.sha}  {\rm H}^1_\fl(K,A_K)\stackrel{\sim}{\longrightarrow} {\rm Ext}(\GGr(A^{\prime 0}), \Q/\Z)={\rm Hom}(\pi_1(\GGr(A^{\prime 0})), \Q/\Z)  ,\end{equation}
where $A^{\prime 0}$ is the identity component of $A^\prime$, $\GGr$ denotes the perfection of the Greenberg realization functor and $\pi_1(\GGr(A^{\prime 0}))$ is the fundamental group of the proalgebraic groups $\GGr(A^{\prime 0})$ (see  \ref{sec.proalg}).
The existence of the pairing was proved by Shafavarevich for the prime-to-$p$ parts, by  B\'egueri in \cite{Beg} for $K$ of characteristic $0$ and by Bester and the author, respectively in \cite{Bes} and \cite{Ber}, for the equal positive characteristic case.

In the first section of the paper we recall the original construction of Shafarevich's duality due to B\'egueri (for $K$ of mixed characteristic). In the second section we slightly modify B\'egueri's construction using rigidificators. 
The latter construction works in any characteristic.
In the third section we construct a morphism as in \eqref{eq.sha} via the relative Picard functor. We  show that it always coincides with  the modified B\'egueri construction and hence with Shafarevich duality for $K$ of characteristic $0$. In the characteristic $p$ case it coincides with Shafarevich duality on the prime-to-$p$ parts. The analogous result for the $p$ parts, although expected, is still open.

%%%%%%%%%%%%%
\section{Shafarevich's duality}\label{sec.sha}
%%%%%%%%%%%

\subsection{Proalgebraic groups and Greenberg realization}\label{sec.proalg}
%%%%%%%

For the theory of proalgebraic groups over the algebraically closed field $k$ we refer to  \cite{SePRO} and \cite{Oo}, II.7.  For the Greenberg realization functor $\Gr(-)$ we refer to \cite{BLR}, 9.6. For its perfection $\GGr(-)$ see \cite{Mi} III, \S 4.  The functor $\Gr(-)$ associates to any smooth group scheme $Y$ of finite type over $R$ a proalgebraic group  scheme $\Gr(Y)$ over $k$ (in the sense of \cite{Oo}) and its perfection  $\GGr(Y)$ is a proalgebraic group  in the sense of \cite{SePRO}.

In the category of proalgebraic groups, the component group functor $\pi_0$ admits a left derived functor $\pi_1$ which is left exact. 
We list below  well-known facts  used in this paper.  Simply assuming them,  the reader will be able to follow the proofs even if he/she is not familiar with the theory of proalgebraic groups. 

\begin{itemize}
\item[i)] If $Y$ is a smooth group scheme of finite type, $\pi_0(\GGr(Y)$ will coincide with the component group of the special fibre of $Y$ and $\pi_1(\GGr(Y))$ is a profinite group.
\item[ii)] If $P$ is a proalgebraic group and $P^0$ is its identity component, then $\pi_1(P)=\pi_1(P^0)$.
\item [iii)] A short exact sequence of smooth $R$-group schemes of finite type $0\to Y_1\to Y_2\to Y_3\to 0$  provides a long exact sequence of profinite groups (cf. \cite{SePRO}, 10.2/1) 
\[0\to \pi_1(\GGr( Y_1))\to \pi_1(\GGr/(Y_2))\to \pi_1(\GGr( Y_3))\to  \pi_0(\GGr (Y_1))\to \pi_0(\GGr (Y_2))\to \pi_0(\GGr( Y_3))\to 0,\]
\end{itemize}

\subsection{The component group of a torus}\label{sec.comptori}
Recall that the N\'eron model $T$ of a torus $T_K$ is locally of finite type, but, in general, not of finite type over $R$. Let $X_K$ denotes the group of characters of $T_K$. Then $T$ is of finite type, i.e., its component group is torsion, if $X_K(K)=0$ (cf. \cite{BLR}, 10.2/1).

\begin{lemma}\label{lem.comptori}
Let  $f\colon T_{1,K}\to T_{2,K}$  be an isogeny of tori with  kernel a finite group scheme $F_K$. The group ${\rm H}^1 _\fl(K, F_K)= T_{2,K}(K)/T_{1,K}(K)$ has a canonical proalgebraic structure.
\end{lemma}
\proof  (Cf.  \cite{Beg}, 4.3.)
Let $X_{i,K}$ be the character group of $ T_{i,K}$ and $T_i$ the N\'eron model, $i=1,2$. 
Denote by $T_{1,K}^{(d)}$ the torus whose group of characters is the constant free group $X_{1,K}(K)$. Similar for  $T_{2,K}^{(d)}$. They are split tori with the same component group, say  $\Z^r$. Furthermore the isogeny $f$ induces an isogeny
$f^{(d)}\colon T_{1,K}^{(d)}\to T_{2,K}^{(d)} $ that is injective on component groups.
The torus $T_{i,K}'$,  kernel of the quotient map $ T_{i,K}\to T_{i,K}^{(d)}$, admits a N\'eron model of finite type because its group of characters is $X_{i,K}'=X_{i,K}/ X_{i,K}(K)$. Hence, using the exact sequences $\pi_0( T_i')\to \pi_0(T_i)\to \pi_0(T_i^{(d)})\to 0$, one sees that the kernel and the cokernel of the homomorphism $ \pi_0(T_1)\to \pi_0(T_2)$ are finite groups.
The identity components of the N\'eron models $T_i$ are smooth group schemes of finite type (\cite{SGA3}, VI A, 2.4). Hence their perfect Greenberg realizations are proalgebraic groups. Let us denote by $P$ the cokernel of the map $\GGr(T_1^0)\to \GGr(T_2^0)$. Now, the cokernel  of the map $\GGr(T_1)\to \GGr(T_2)$ is extension of the
finite group $\pi_0(T_2)/ \pi_0(T_1)$ by  the quotient of $P$ by a finite constant group. Hence it is proalgebraic.  
\qed

We will use in the next sections the following result: 

\begin{lemma}\label{lem.torustor}
Let $0\to T\to G\to A^0\to 0$ be an exact sequence of smooth $R$-group schemes over $R$ where $T$ is the N\'eron model of a torus and $A^0$ is the identity component of the N\'eron model of an abelian variety. It induces a morphism $\pi_1(\GGr (A^0)) \to \pi_0(T)_\tor$ where the index $\tor$ indicates the torsion subgroup. 
\end{lemma}
\proof
If $T$ is of finite type, the proof is immediate. Indeed $\pi_0(T)_\tor=\pi_0(T)=\pi_0(\GGr(T))$ is finite and, on applying the Greenberg functor to the sequence above, we get an extension
\[
0\to \GGr(T)\to \GGr(G)\to \GGr(A^0)\to 0.
\]
The desired map follows then from long exact sequence of the $\pi_i$'s.

Suppose that $T$ is  locally of finite type. As $k$ is algebraically closed, $\pi_0(T)=\pi_0(T)_\tor\times \pi_0(T)_\fr$.  Let $T^\ft$ be the maximal subgroup of $T$ such that its component group is finite, i.e., $T^\ft$ contains the identity component $T^0$ and  $\pi_0(T^\ft)=\pi_0(T)_\tor$. The extension $G$ is the push-out  along the inclusion $T^\ft\to T$ of a unique (up to isomorphism) extension  $0\to T^\ft\to G^\ft\to A^0\to 0$  because  the quotient $T/T^\ft$  satisfies the hypothesis in \cite{SGA7}, \S5.7, 5.5. 
Hence we can proceed as done above with $G^\ft$ in place of $G$ obtaining a map  $\pi_1(\GGr (A)) \to \pi_0(\GGr(T^\ft))=\pi_0(T)_\tor $ .
\qed

\subsection{B\'egueri's construction}
Suppose in this section that $K$ has characteristic $0$.
Given $K$-schemes $Z_K,S_K$ let us denote by $Z_{S_K}$ the fibre product $Z_K\times_K S_K$, viewed as a scheme over $S_K$.

Let $X_K$ be a  $K$-torsor under $A_K$ and let $n$ be a  positive integer such that $n\cdot X_K$ is trivial; the order of $X_K$ is the minimum among such integers.  The torsor $X_K$ corresponds to a $n$-torsion 
element in  ${\rm H}^1_\fl(K,A_K)={\rm Ext}^1(\Z,A_K)$ 
and hence to an extension 
\begin{equation}\label{eq.etaB} 0\to A_K\to B_K\to \Z\to 0,
\end{equation}
that is the pull-back along $\Z\to \Z/n\Z$ of a, not unique, extension
\begin{equation}\label{eq.eta}\eta\colon 0\to A_K\stackrel{\alpha}{\to} E_K\to \Z/n\Z\to 0,\end{equation}
precisely $X_K$ is the fibre at $1\in  \Z/d\Z$.  Let us denote by 
\begin{equation}\label{eq.etan}
\eta_n  \colon 0\to{}_n A_K\to {}_n E_K\to \Z/n\Z\to 0
\end{equation}
the sequence of $n$-torsion subgroups.

Consider the exact sequence 
\begin{equation}\label{eq.muvE}
0\to {\mu_n}\to V_{{}_nE_K}^*\to \underline{\rm Ext}^1(E_K,\G_m)_{{}_nE_K}\stackrel{\tau}{\to}
  \underline{\rm Ext}^1(E_K,\G_m)=A^\prime_K \to 0
\end{equation}
(cf. \cite{Beg}, 2.3.2)
where we denote by $ V_{{}_nE_K}^*$ the Weil restriction torus (\cite{BLR}, 7.6) \[\Re_{{}_nE_K/K}(\G_{m, {}_nE_K})=\underline{\rm Mor}({}_nE_K,\G_{m,K})\]
representing the functor that associates to a $K$-scheme $S$ the group $\G_{m,K}(S\times_K  {}_nE_K  )$. Observe that  \[\mu_n=\uline{\rm Hom}(\Z/d\Z,\G_m)=\uline{\rm Hom}(E_K,\G_m).\]
The first map in \eqref{eq.muvE} maps a homomorphism $f\colon E_K\to \G_{m,K}$ to its restriction to ${}_nE_K$,  the second arrow sends $g\in \G_m( {}_nE_K ) $ to  (the isomorphism class of)  the trivial extension endowed with the section $g$, the latter map forgets the rigidification along ${}_nE_K$.

We describe now B\'egueri's construction of Shafarevich's duality following \cite{Beg}. 
Let $F_K$ be a finite $K$-group scheme and $F_K^D$ its Cartier dual. There is  
 a short exact sequence (cf. \cite{Beg}, 2.2.1)
\begin{equation}\label{eq.muZ2}
0\to F_K^D\to V_{F_K}^*\to \underline{\rm Ext}^1(F_K,\G_m)_{F_K}\to 0,
\end{equation}
where the first map forgets the group structure and the second map associates to a $f\in \G_{m,K}(F_K)$ the trivial extension endowed with the rigidification induced by $f$.
Finally we introduce  the following exact sequence (cf. \cite{Beg}, 2.3.1)
\begin{equation}\label{eq.ExtAAprime}
0\to V_{{}_nA_K}^*\to \underline{\rm Ext}^1(A_K,\G_m)_{{}_nA_K}\stackrel{\tau}{\to}
 A^\prime_K\to 0.
\end{equation}
At first B\'egueri constructs in \cite{Beg}, 8.2.2, a map \begin{equation}\label{eq.Gamma}\Gamma\colon {\rm H}^1_\fl(K,{}_nA_K)\to {\rm Ext}^1(\GGr(A'),{\uuline {\rm H}}^1(K,\mu_n) )\end{equation} as follows:  any element in ${\rm H}^1_\fl(K,{}_nA_K)$ corresponds to a sequence $\eta_n$ as in \eqref{eq.etan}. Consider  now the diagram
\begin{equation}\label{eq.diaBeg}
\xymatrix{
\mu_n\ar@{=}[r]\ar@{^{(}->}[d]&\mu_n\ar@{^{(}->}[d]&0\ar[d]\\
V_{\Z/n\Z}^*\ar[r]\ar[d]^{v_1}&V_{{}_nE_K}^*\ar[r]\ar[d]^{v_2}&V_{{}_nA_K}^*\ar[d]^{v_3}\\
 \underline{\rm Ext}^1(\Z/n\Z,\G_m)_{\Z/n\Z}\ar[d]\ar[r]& \underline{\rm Ext}^1(E_K,\G_m)_{{}_nE_K}\ar@{->>}[d]^\tau\ar[r]& \underline{\rm Ext}^1(A_K,\G_m)_{{}_nA_K}\ar@{->>}[d]\\
0& A'_K\ar@{=}[r]& A'_K&}
\end{equation}
where the rows are complexes and the vertical sequences are those in \eqref{eq.muZ2}, for $F_K=\Z/n\Z$,
 \eqref{eq.muvE}, \eqref{eq.ExtAAprime}, respectively. Since $K$ has characteristic $0$, the second row consists of tori, the third row consists of semi-abelian varieties. Hence they all admit N\'eron models. On passing to the perfection of the Greenberg realization of the N\'eron models and considering the cokernels of the maps induced by $v_1,v_2,v_3$, one gets a complex of proalgebraic groups\footnote{Note that one applies here Lemma~\ref{lem.comptori} to the isogeny $v_1$.}
\begin{equation}\label{eq.ZEA}
 0\to  \uuline{\rm H}^1(K,\mu_n)   \to  \uuline{\rm Ext}^1(E_K,\G_m)\to \GGr(A^\prime)\to 0; \end{equation}
it is indeed an exact sequence because on $k$-points it induces the exact sequence
\[ 0\to{ \rm H}^1_\fl(K,\mu_n)={\rm Ext}^1(\Z/d\Z,\G_m)\to {\rm Ext}^1(E_K,\G_m)\to A'(R)= {\rm Ext}^1(A_K,\G_m)\to 0. \]
We have then associated to \eqref{eq.etan} an extension of $\GGr(A')$ by ${\uuline{\rm H}}^1(K,\mu_n)$: this is the image of \eqref{eq.etan} via  $\Gamma$.

The homomorphism 
\begin{equation}\label{eq.psid}
\psi_n\colon {\rm H}^1_\fl(K,{}_nA_K)\longrightarrow {\rm Ext}^1(\GGr(A^{\prime 0}),\Z/n\Z)
 \end{equation}
in \cite{Beg}, 8.2.3, is then obtained by  applying  first $\Gamma$, then the pull-back along
$\GGr(A^{\prime 0})\to \GGr(A^{\prime})$ and, finally, the push-out along ${\uuline{\rm H}}^1(K,\mu_n)\to \pi_0({\uuline{\rm H}}^1(K,\mu_n)  ) =\Z/n\Z. $

Let us denote by
\begin{equation}\label{eq.shaX}
\psi_n(\eta_n)\colon \quad  0\to \Z/n\Z \to W(X_K)\to \GGr(A^{\prime 0})\to 0  \end{equation}
the image of \eqref{eq.etan} via $\psi_n$, 

Recall now that (cf. \cite{SePRO}, 5.4)
\[{\rm Ext}(\GGr(A^{\prime 0}), \Q/\Z)={\rm Hom}(\pi_1(\GGr(A^{\prime 0})), \Q/\Z)={\rm Hom}(\pi_1(\GGr(A^\prime)), \Q/\Z) .\]
In terms of homomorphisms of profinite groups, the extension \eqref{eq.shaX} corresponds to a map
 \begin{equation}\label{eq.utau}
u^\tau=u^\tau_{X_K} \colon  \pi_1(\GGr(A^{\prime 0}))\longrightarrow  \pi_0(\uuline{\rm H}^1(K,\mu_n)  )=\Z/n\Z\subset \Q/\Z
\end{equation}
deduced from \eqref{eq.ZEA} (or the same, from the pull-back of 
\eqref{eq.ZEA} along  $\GGr(A^{\prime 0})\to \GGr(A')$)
 via the long exact sequence of $\pi_i$'s.

It is proved in \cite{Beg}, 8.2.3 that:  a)\ the extension $\psi_n(\eta_n)$ in \eqref{eq.shaX}
 depends only on the sequence \eqref{eq.etaB}, i.e., on the torsor $X_K$,  b)\ its formation behaves well w.r.t. the inclusions $ {}_n{\rm H}^1_\fl (K,{}_nA_K)\to    {}_{n'}{\rm H}^1_\fl( K,{}_{n'}A_K)$, and $\Z/n\Z\to \Z/n'\Z\to \Q/\Z$ for $n|n'$, and c)\  it provides the duality in \eqref{eq.sha}.  Hence:
\medskip

\begin{verse}
\fbox{
Shafarevich's duality in  \eqref{eq.sha} maps the torsor $X_K$ to the homomorphism $ u^\tau_{X_K}$ in \eqref{eq.utau}.}
\end{verse}

\subsubsection{An alternative construction of  $\psi_n(\eta_n)$ in \eqref{eq.shaX} (in view of  further applications).}

The kernel of $\tau$  in \eqref{eq.ExtAAprime} is a torus, that for brevity we denote by $T^\tau_K$. Let $T^\tau$ be  its N\'eron model. We have an exact sequence 
\begin{equation*}\label{seq.torustau}
0\to T_K^\tau\to \underline{\rm Ext}^1(E_K,\G_m)_{{}_nE_K}\stackrel{\tau_K}{\to}
 A^\prime_K\to 0.
 \end{equation*}
which extends to an exact sequence of  N\'eron models 
\begin{equation}\label{seq.torustauN}
0\to T^\tau\to j_*\underline{\rm Ext}^1(E_K,\G_m)_{{}_nE_K}\stackrel{\tau}{\to}
 A^\prime \to 0.
 \end{equation}
On applying the perfection of the Greenberg functor we get an exact sequence
\begin{equation}\label{eq.Grj}
0\to\GGr( T^\tau)\to \GGr(j_* \underline{\rm Ext}^1(E_K,\G_m)_{{}_nE_K})\stackrel{\tau}{\to}
\GGr( A^\prime)\to 0
\end{equation}
where the first two groups are not proalgebraic,  in general. Nevertheless,  on applying the  Greenberg functor to the morphism of
N\'eron models $j_*V_{{}_nE_K}^*\to  T^\tau$ we get a homomorphism whose cokernel  is a proalgebraic group (cf. Lemma~\ref{lem.comptori}) and 
whose group of $k$-points is ${\rm H}^1(K,\mu_n)$;  we will write
\begin{equation}\label{eq.GrH}
\GGr( j_*V_{{}_nE_K}^*  )\to \GGr(T^\tau)\stackrel{h^\tau}{\to}\uuline{\rm H}^1(K,\mu_n) \to 0.
\end{equation}

Now  take  the push-out  of \eqref{eq.Grj}
 along $h^\tau$;  by construction, the resulting exact sequence
 is the one in \eqref{eq.ZEA}, i.e., the image of 	\eqref{eq.etan} via $\Gamma$. 
Hence, if one considers the pull-back of \eqref{eq.Grj} along $\GGr(A^{\prime 0})\to \GGr(A^\prime)$,
\begin{equation}\label{eq.Grj0}
0\to\GGr( T^\tau)\to  U {\to}
\GGr( A^{ \prime 0})\to 0,
\end{equation}
and then the push-out of \eqref{eq.Grj0} along the composition of maps
\[  \GGr(T^\tau)\stackrel{h^\tau}{\to}\uuline{\rm H}^1(K,\mu_n)\to\pi_0( \uuline{\rm H}^1(K,\mu_n) )=\Z/n\Z , \]
one gets the extension $\psi_n(\eta_n)$ in \eqref{eq.shaX}, i.e., the image of $X_K$ via Shafarevich's duality

Thanks to this new description of Shafarevich's map, we can  characterize of the map $u^\tau$ in \eqref{eq.utau} as follows:

First consider the pull-back of  \eqref{seq.torustauN} along $A^{\prime 0}\to A'$. As we have seen in the proof of Lemma~\ref{lem.torustor}, this extension is the push-out of an extension 
\[0\to T^{\tau, \ft}\to G\to A^{\prime 0}\to 0\]
where $T^{\tau, \ft}$ is the maximal subgroup scheme  of finite type of $T^\tau$. The pull-back of the sequence in \eqref{eq.Grj}
along $\GGr( A^{\prime 0})\to \GGr(A^\prime)$
is then isomorphic to the push-out of
\begin{equation}\label{eq.zetataup}
0\to\GGr( T^{\tau, \ft})\to \GGr(G)\to \GGr(A^{\prime 0})\to 0\end{equation}
along the composition of maps $h^{\tau, \ft}\colon \GGr( T^{\tau, \ft})\to \GGr( T^\tau)\stackrel{h^\tau}{\to}  \uuline{\rm H}^1(K,\mu_n)$. Hence 
\begin{equation}\label{eq.w}
\fbox{$u^\tau_{X_K} =\pi_0(h^{\tau, \ft})\circ w$}
\end{equation}
where the homomorphism $w\colon \pi_1(\GGr(A^\prime))\to  \pi_0(\GGr(T^{\tau, \ft}))$ is deduced from 
\eqref{eq.zetataup}  via the long exact sequence of $\pi_i$'s. 

%%%%%%
\section{An alternative construction using  rigidificators}\label{sec.rig}
%%%%%%%%%%%%
Let $X_K$ be a torsor under an abelian variety $A_K$.
We will see in this section how  the homomorphism $u^\tau$ in \eqref{eq.utau}  (and in \eqref{eq.w}) can be constructed using a rigidificator $x_K$ of the relative Picard functor $\Pic_{X_K/K}$.

\subsection{Rigidificators}
A rigidificator $x_K$ of the relative Picard functor $\Pic_{X_K/K}$  is  a finite closed subscheme $x_K$ of $X_K$ such that for any $K$-scheme $S_K$ the map 
\[\Gamma(X_{S_K},{\mathcal O}_{X_{S_K}})\to \Gamma(x_{S_K},{\mathcal O}_{x_K\times_K{S_K}})  \]
is injective. (cf. \cite{Ray}, 2.1.1).

\begin{lemma}\label{lem.index-period} Let $X_K$ be a torsor under $A_K$ of order $d$. Let $d'$ be the separable index of $X_K$, i.e., the greatest common divisor of the degrees of its finite separable splitting extensions. Then  $d|d'$ and they have the same prime factors. 
If $A_K$ is an elliptic curve, then $d=d'$ and the index is indeed the degree of a minimal separable splitting extension.
\end{lemma}
\proof One can invoke  \cite{LT}, Proposition 5. Alternatively one proves, using the restriction and corestriction maps, that $n\cdot X_K=0$  if  $X_K$ becomes trivial over a finite separable extension $K'/K$ of degree $n$. Hence $d|n$. Suppose now given separable extensions $K\subseteq L\subseteq L'$ with $(d,[L':L])=1$ and $X_{L'}=0$. Then $[L':L]\cdot X_{L}=0$ in ${\rm H}^1_\fl(L,A_{L})$. However the order of $X_{L}$ in ${\rm H}^1_\fl(L,A_{L})$ divides $d$ and  thus $X_{L}=0$. Hence $d,d'$ have the same prime factors.

For the latter assertion on elliptic curves see \cite{LT}, p. 670, or  \cite{Licht}, Theorems 1 \& 4, or \cite{To}, 2.1.2.
\qed

\begin{remark}\label{rem.tori} Let $x_K=\spec(K')$ with $K'/K$ a finite separable extension of degree $n$. Then the torus $V_{x_K}^*$ is the Weil restriction $\Re_{K'/K}(\G_{m,K'})$, it has component group isomorphic to $\Z$ and the closed immersion  $\G_{m,K}\to  V_{x_K}^*$ (i.e., the inclusion $K^*\subset K^{\prime *}$ on $K$-sections) induces the $n$-multiplication $n\colon \Z\to \Z$ on component groups of N\'eron models. 
\end{remark}

\subsection{The alternative construction of Shafarevich's duality}
The main idea is to use a rigidificator $x_K$ in place of ${}_nE_K$ in \eqref{eq.muvE}. 
The advantage is that the new construction works even for $K$ of positive characteristic; in this case we choose $x_K$  \'etale so  that  $V_{x_K}^*=\Re_{x_K/K}(\G_{m, x_K})$ is still a torus.

Observe that a rigidificator $x_K$ is a closed subscheme of $E_K$ and the homomorphism
 \begin{equation}\label{eq.muvclosed}
\mu_n=\underline{\rm Hom}(E_K,\G_m)\to V_{x_K}^*\end{equation}
is still a closed immersion.  
Indeed any homomorphism $f\colon E_K\to \G_m$ factors through $\rho\colon E_K\to \Z/d\Z$ and if $f_{|x_K}=0$ then $f_{|X_K}=0$ because $x_K$ is a 
rigidificator. However $X_K$ is the fibre at $1$ of $\rho$ and hence also $f=0$. 
We then have  an exact sequence 
\begin{equation}\label{eq.muvEx}
0\to {\mu_n}=\underline{\rm Hom}(E_K,\G_m) \to V_{x_K}^*\to \underline{\rm Ext}^1(E_K,\G_m)_{x_K}\stackrel{ }{\to}
 A^\prime_K\to 0 .
\end{equation}

More generally, we will say that a finite \'etale subscheme $Z_K$ of $E_K$ satisfies property $(*)$ if
\begin{verse}
 $ (*)$\quad   the canonical map
$\mu_n=\underline{\rm Hom}(E_K,\G_m)\to V_{Z_K}^*$  is a closed immersion.
\end{verse}
For any such \'etale subscheme $x_K$ we can construct an exact sequence as in \eqref{eq.muvEx}.

Denote by $T_K^x$ the torus $V_{x_K}^*/\mu_n $ and omit the exponent $x$ if the rigidificator $x_K$ is fixed. 
The sequence \eqref{eq.muvEx} induces an exact sequence 
\begin{equation}\label{eq.TExtA}
0\to T_K\to \underline{\rm Ext}^1(E_K,\G_m)_{x_K}\stackrel{ }{\to}
 A^\prime_K\to 0, 
\end{equation} 
and hence a n exact sequence (cf. Lemma~\ref{lem.torustor}/proof)
\begin{equation}\label{eq.TExtAN}
0\to T^\ft\to G_1\stackrel{ }{\to}
 A^{\prime 0}\to 0, 
\end{equation} 
where $T^\ft$ is the maximal subgroup of finite type of the N\'eron model $T$ of $T_K$.
Consider now the cokernel 
 \begin{equation}\label{eq.hproalg}\GGr(j_*V_{x_K}^*)\stackrel{g^x}{\to}\GGr(T)\stackrel{h}{\to}  \uuline{\rm H}^1(K,\mu_n)\to 0\end{equation}
of the homomorphism between  the Greenberg realizations of the  N\'eron models of $V_{x_K}^*$ and $T_K$; by Lemma~\ref{lem.comptori} it 
is a proalgebraic group whose group of $k$-points is  $ {\rm H}^1_\fl(K,\mu_n)$.

\begin{lemma}
The proalgebraic group $\uuline {\rm H}^1(K,\mu_n)$ in \eqref{eq.hproalg} does not depend on the \'etale finite subscheme $x_K$ chosen to construct it. In particular it coincides with the one in \eqref{eq.muvE}.\end{lemma}
\proof
Let $x_K\subset y_K$ be finite \'etale subschemes of $E_K$ satisfying $(*)$. 
We have canonical morphisms  $    f^{ V}\colon  V_{y_K}^*\to      V_{x_K}^* $,    $f^{ T}\colon T_K^y\to T_K^x$, such that $f^{ T}\circ g^y=g^x\circ f^{\rm V}$. Hence the proalgebraic group constructed in \eqref{eq.hproalg} for $x_K$  is canonically isomorphic to the one constructed via $y_K$. 
\qed

In order to provide a more useful description of the map $u^\tau$ in \eqref{eq.utau}, consider  the perfect Greenberg realization of  \eqref{eq.TExtAN}
\begin{equation}\label{eq.TExtANG}
0\to \GGr(T^\ft)\to \GGr(G_1)\stackrel{ }{\to}
\GGr( A^{\prime 0})\to 0, 
\end{equation} 
and then its push-out  along the composition of maps
\begin{equation}\label{eq.hprime}h^\ft\colon \GGr(T^\ft)\to \GGr(T)\stackrel{h}{\to} \uuline{\rm H}^1(K,\mu_n).\end{equation}
We obtain an exact sequence 
\begin{equation}\label{eq.zeta}
\zeta\colon \quad 0\to  \uuline{\rm H}^1(K,\mu_n)\to W' \to \GGr(
 A^{\prime 0})\to 0 
\end{equation} 
and hence a homomorphism
\[u_{X_K}=u\colon \pi_1(\GGr(A^\prime))\to \pi_0(  \uuline{\rm H}^1(K,\mu_n)   )=\Z/n\Z\subset \Q/\Z\]
such that 
\begin{equation}\label{eq.u} u=\pi_0(h^\ft)\circ u^\ft,\end{equation} 
where $ u^\ft\colon \pi_1(\GGr(A^\prime))\to \pi_0(\GGr(T^\ft))=\pi_0(T)_\tor$
is deduced from the long exact sequence of $\pi_i$'s of \eqref{eq.TExtANG}.

\begin{proposition}\label{pro.sha0}
The association $X_K\mapsto u_{X_K}$ provides a homomorphism
\[{\rm H}^1(K, A_K)\to {\rm Hom}(\pi_1(\GGr(A^\prime)),\Q/\Z).   \]  
If $char(K)=0$ it is Shafarevich duality in \eqref{eq.sha}, i.e., $u_{X_K}=u^\tau_{X_K}$. 
\end{proposition}
\proof
We start showing that the construction of $u\colon \pi_1(\GGr(A'))\to \Q/\Z$ in  \eqref{eq.u} does not depend on the choices of  $x_K$, $n$ and  $\eta\in {\rm Ext}^1(\Z/n\Z,A_K)$ above $X_K$.

At first we see that $u$ does not depend on  the \'etale finite closed subscheme  $x_K$ of $E_K$ satisfying $(*)$. 
Let $x_K\subset y_K$ be two \'etale  subschemes of $E_K$ satisfying $(*)$. 
Denote by $T_K^x, h^x, h^{\ft, x}, u^{\prime x}, u^x$ respectively the torus  in \eqref{eq.TExtA}, the maps in \eqref{eq.hproalg},  \eqref{eq.hprime} and \eqref{eq.u} for $x_K$, and similarly for $y_K$.
We have a canonical morphism of tori $ T_K^y\to T_K^x$ and it induces a morphism $\beta\colon T_K^{y,\ft}\to  T_K^{x,\ft }$ between the maximal subgroups of finite type. Denote by $\beta^\prime\colon \GGr(T^{y, \ft})\to \GGr(T^{x,\ft})$ the corresponding map on Greenberg realizations of N\'eron models. It holds $\beta^\prime  \circ h^{x,\ft}= h^{y,\ft}$.  One has  $\pi_0(h^{y,\ft})=\pi_0(h^{x,\ft})\circ \pi_0(\beta')$.
Furthermore the sequence  \eqref{eq.TExtANG}  for $x_K$ is the push-out along $\beta^\prime$ of the sequence \eqref{eq.TExtANG} for $y_K$. Hence $u^{x,\ft}=\pi_0(\beta^\prime )\circ u^{y,\ft}$. We conclude then that
\begin{equation}\label{eq.uxuy}  u^x= \pi_0(h^{x,\ft}) \circ u^{x,\ft} = \pi_0(h^{x,\ft}) \circ \pi_0(\beta^\prime )\circ u^{y,\ft}=\pi_0(h^{y,\ft}) \circ u^{y,\ft}=u^y.
\end{equation}

Let now $n,\hat n$ be positive integers such that $n\cdot X_K=0$ and $n|\hat n$. We can consider the pull-back $\hat \eta$ of $\eta$ in \eqref{eq.eta} along the projection $\Z/\hat n\Z\to \Z/n\Z$. If we proceed with $\hat \eta$ as we have done for $\eta$, we get a   map  $\hat u\colon \pi_1(\GGr(A^\prime))\to \Q/\Z$. Observe that the $2$-fold extension \eqref{eq.muvEx} for $\hat\eta$ is the push-out along $\mu_n\to \mu_{\hat n}$ of  \eqref{eq.muvEx} and that the map
$ \pi_0  (  \uuline{\rm H}^1(K,\mu_n))\to  \pi_0  (  \uuline{\rm H}^1(K,\mu_{\hat n}))$  is the inclusion $\Z/n\Z\to \Z/\hat n\Z$.  It is now immediate to check that the maps  $\hat u$ and $u$ coincide. 

We have then obtained a map
\begin{equation}\label{eq.etau}{\rm Ext}^1(\Z/n\Z, A_K)\to{\rm Hom}(\pi_1(\GGr(A')),\Q/\Z),\qquad \eta\mapsto u. \end{equation}
To check that it is indeed an homomorphism, observe that it
 is functorial in $A_K$. Furthermore we could repeat the construction with any finite constant group $F_K$ in place of $\Z/n\Z$ obtaining in this way a map
 \[{\rm Ext}^1(F_K,A_K)\to {\rm Hom}( \pi_1(\GGr(A^\prime)),\pi_0(\uuline{\rm H}^1(K,F_K^D)) )
\]
with $F_K^D$ the Cartier dual of $F_K$. This construction is functorial in $F_K$. 
The functoriality results is sufficient  to conclude that the map in \eqref{eq.etau} is a homomorphism, because the Baer's sum of two extensions as in \eqref{eq.eta}
is done first by  taking the direct sum of the two extensions, then by applying the push-out along the multiplication of $A_K'$ and finally by applying the pull-back along the diagonal $\Z/n\Z\to \Z/n\Z\oplus \Z/n\Z$.

Suppose now $n$ and $x_K$  fixed. We show now that the  map $u$ is trivial if $X_K$ is trivial, i.e.,  the map in \eqref{eq.etau} factors through ${\rm H}^1_\fl(K,A_K)$.  Suppose $X_K$ is trivial and choose a $K$-point $x_K$ of $X_K$. In particular, $V_{x_K}^*=\G_{m,K}$, $T_K=\G_{m,K}$ and  $\pi_0(T)=\Z$. 
 Hence $T^\ft=\G_{m,R}$, the homomorphism $u^\ft\colon \pi_1( \Gr(A'))\to \pi_0(\GGr(T^\ft))=0$ is the zero map and  $u=0$.

Suppose now  $char(K)=0$. To see that the homomorphism  $X_K\mapsto u_{X_K}$ is Shafarevich's duality, it is sufficient to check that 
the homomorphisms   $u^\tau$ in \eqref{eq.utau}
and $u$ in \eqref{eq.u} coincide. 
Consider then  a  finite separable extension  $K'/K$ splitting \eqref{eq.etan}   and  a point $x_K=\spec(K')$ of ${}_n E_K$  above $1$. It is  a rigidificator of $\Pic_{X_K/K}$.  Pose $y_K={}_n E_K$. 
Then, using notations as above,    $u^y$ coincides with the map $u^\tau$ in \eqref{eq.w} and one can repeat the arguments in \eqref{eq.uxuy}.
\qed

\begin{remark}
The original construction by B\'egueri works only for $K$ of characteristic zero because in the case of positive characteristic the scheme $  V_{{}_n E_K}^*$ (and hence $T^\tau_K$) might not be a torus; in particular it might not admit a N\'eron model. 
The  construction via rigidificators described in this section works in any characteristic. For $char(K)=p$ it is not clear that it provides Shafarevich duality. We will see in Proposition \ref{pro.shap} that this is the case on the prime-to-$p$ parts. \end{remark}

%%%%%%
\section{A construction via the Picard functor}
%%%%%

In this section we present a third possible construction of a homomorphism as in \eqref{eq.sha} which makes use the relative Picard functor. 
We will see that it always coincide with the one in Proposition~\ref{pro.sha0} and hence with Shafarevich's duality in the characteristic $0$ case.

Let $X_K$ be a torsor under $A_K$ and $x_K=\spec(K')$ a rigidificator of $\Pic_{X_K/K}$ with $K'/K$ a finite separable extension.  It exists by Lemma~\ref{lem.index-period}. No assumption on the characteristic of $K$ is made.

Consider the usual exact sequence (cf. \cite{Ray}, 2.4.1)
 \begin{equation}\label{eq.GVPic}
0\to  V_{X_K}^*\to V_{x_K}^*\to ({\rm Pic}_{X_K/K},x_K)^0\to
 A^\prime_K\to 0 .
\end{equation}
Observe that  $ V_{X_K}^*=\Re_{X_K/K}(\G_{m,X_K})=\G_{m,K}$ (\cite{Ray}, 2.4.3),  $ V_{x_K}^*$ is a torus and hence  so too is $N_K:=V_{x_K}^*/\G_{m,K}$. Denote by $N$ its N\'eron model.  Observe that it follows from Remark~\ref{rem.tori} that the component group of  $N$ is cyclic of order $n$, hence the perfection of its Greenberg realization is a proalgebraic group.

We proceed as in the previous section, first by passing to  N\'eron models and then  applying the Greenberg  realization to the sequence 
\begin{equation}\label{eq.N}
0\to N_K\to  (\Pic_{X_K/K},x_K)^0 \stackrel{h_K}{\to}  A^\prime_K   \to 0
\end{equation}
so that we obtain an exact sequence of proalgebraic groups
\begin{equation}\label{eq.GrN}
0\to \GGr(N)\to \GGr(j_* (\Pic_{X_K/K},x_K)^0 )\stackrel{h}{\to}  \GGr(A^\prime)   \to 0
\end{equation}
and hence a  homomorphism 
\begin{equation}\label{eq.v}
v=v_{X_K}\colon \pi_1(\GGr(A^\prime))\longrightarrow  \pi_0(\GGr(N))=\Z/n\Z.
\end{equation}

In order to compare this construction with the (modified) B\'egueri construction of the previous section, i.e., in order to compare the maps $u$ \eqref{eq.u} and $v$ in \eqref{eq.v},  we consider the following diagram
\begin{equation}\label{eq.diagram}
\xymatrix{
0\ar[r]& [V_{x_K}^*/\mu_n]=T_K\ar[r]\ar[d]^{t_K}&  \underline{\rm Ext}^1(E_K,\G_m)_{x_K}\ar[d]^{f_K}
                         \ar[r]& A^\prime_K  \ar@{=}[d]  \ar[r]& 0\\
0\ar[r]& [V_{x_K}^*/ V_{X_K}^*]= N_K\ar[r]& (\Pic_{X_K/K},x_K)^0 \ar[r]^(0.6){h_K}&  A^\prime_K   \ar[r]& 0}
\end{equation} 
where the upper sequence is \eqref{eq.TExtA}, the lower one is \eqref{eq.N}
and  $f_K$ associates to a $\G_m$-extension $L_K$ of $E_K$ endowed with a 
$x_K$-section $\sigma$  its restriction (as torsor) to $X_K$ endowed with the trivialization along $x_K$
induced by $\sigma$.  The morphism $t_K$ is surjective and its kernel is $\G_{m,K}=V_{X_K}^*/\mu_n=\G_{m,K}/\mu_n$.

Consider now the induced diagram on N\'eron models. 
\[
\xymatrix{
0\ar[r]& T^\ft\ar[r]\ar[d]^{}\ar@/_1pc/[dd]_{t^\ft}& G_1\ar[d]^{}
                         \ar[r]& A^{\prime 0}  \ar[d]  \ar[r]& 0\\& T\ar@{^{(}->}[r]\ar[d]& j_* \underline{\rm Ext}^1(E_K,\G_m)_{x_K}\ar[d]^{f}
                         \ar[r]& A^\prime  \ar@{=}[d]  \ar[r]& 0\\
0\ar[r]& N\ar[r]& j_*(\Pic_{X_K/K},x_K)^0 \ar[r]&  A^\prime   \ar[r]& 0}
\]
where the first row is \eqref{eq.TExtAN}.
The homomorphism $u$ in \eqref{eq.u} is the composition of the homomorphism 
$u^\ft\colon \pi_1(\GGr(A^{\prime 0})) \longrightarrow  \pi_0(\GGr(T^\ft))$ (deduced from the upper exact sequence) with the homomorphism 
\[\pi_0(h^\ft)\colon \pi_0(T^\ft)= \pi_0(\GGr(T^\ft))\to \pi_0(\uuline{\rm H}^1(K,\mu_n)   ).\] Now it follows form the above diagram that the map $v\colon \pi_1(\GGr(A^\prime))\longrightarrow  \pi_0(\GGr(N))$ in \eqref{eq.v}, obtained from the lower exact sequence, satisfies
\begin{equation}\label{eq.vuft} v= \pi_0(t^\ft)\circ u^\ft.  
 \end{equation}
In order to explicate the relation between $u$ and $v$, consider the following diagram
\begin{equation*}
\xymatrix{&&&\G_{m,K} \ar@{^{(}->}[d]&\\
0\ar[r]& \mu_n\ar[r]  \ar[d] &V_{x_K}^* \ar[r]\ar@{=}[d]&  T_K \ar[r]\ar[d]^{t_K}&  0\\
0\ar[r]&V_{X_K}^* =\G_{m,K}\ar@{->>}[d]^n\ar[r]& V_{x_K}^*  \ar[r]& N_K\ar[r]& 0\\
& \G_{m,K} &&&}
\end{equation*}
and consider the induced diagram on component groups of N\'eron models
\begin{equation*}
\xymatrix{ & & & \Z\ar@{^{(}->}[d]&\pi_0(T^\ft)=\pi_0(T)_\tor\ar@{_{(}.>}[dl]_\iota
\ar@{.>}[ddl]^{\pi_0(t^\ft)} \\
&&\pi_0(V_{x_K}^*)\ar[r]\ar@{=}[d]&  \pi_0(T) \ar@{->>}[d]_{}&\\
0\ar[r]&\Z  \ar[r]&\pi_0( V_{x_K}^* ) \ar[r]& \pi_0(N)\ar[r]& 0}
\end{equation*}
where $\iota$ is the inclusion map and the vertical  sequence is left exact because $\Z$ is torsion free (cf. \cite{SGA7}, VIII 5.5).

We insert this diagram into a bigger diagram
\begin{equation}\label{eq.diagrambig}
\xymatrix{ 0\ar[r]& \Z\ar@{^{(}->}[d]\ar[r]^n& \Z\ar[r]^{}\ar[d]& \Z/n\Z\ar[r]\ar@{=}[d]&0\\
 0\ar[r]&\pi_0(V_{x_K}^* )\ar[r]\ar@{->>}[d]^{}&  \pi_0(T)\ar[r] ^{q_1~ ~}\ar@{->>}[d]_{q_2}&\pi_0(  \uuline{\rm H}^1(K,\mu_n))  \ar[r]&0\\
&  \pi_0(N) \ar@{=}[r]  & \pi_0(N )& \pi_0(T^\ft)\ar@{.>}[u] _{\pi_0(h^\ft)}\ar@{.>}[l]^{\pi_0(t^\ft)}\ar@{_{(}.>}[ul]_\iota&
}
\end{equation}
where $q_1\circ \iota= \pi_0(h^\ft)$ and $q_2\circ \iota=\pi_0(t^\ft)$.
By Remark~\ref{rem.tori},  $\pi_0(N)=\Z/n\Z$ and the vertical sequence on the left coincides with the upper horizontal sequence. 
Hence the vertical sequence in the middle splits as well as the horizontal sequence in the middle. The identifications $\pi_0(V_{x_K}^* )=\Z$ induces then the identification
\[ \pi_0(N)\cong \Z/n\Z \cong \pi_0 (  \uuline{\rm H}^1(K,\mu_n))\] where the first 
isomorphism maps the image of the class of a uniformizer $\pi'\in K^{\prime *}=V_{x_K}^* (K)$ to the class of $1$, while the second isomorphism maps the class of $1$ to the image of the cohomology class corresponding to a uniformizer $\pi\in K^*=\G_{m,K}(K)$. 

Let $\sigma$ be a section of $q_2$. It holds $q_1\circ \sigma=\id_{\Z/n\Z}$.
Furthermore $\sigma\circ q_2\circ \iota=\iota$ because $\sigma\circ q_2\circ \iota-\iota$ factors through $\Z$ and thus is trivial because $\pi_0(T^\ft)$ is torsion. Hence 
\[\pi_0(h^\fl)=q_1\circ \iota= q_1\circ \sigma \circ q_2\circ \iota=q_2\circ \iota=\pi_0(t^\ft)\]
and hence thanks to  \eqref{eq.vuft} and \eqref{eq.u}, we get
\[v=\pi_0(t^\ft)\circ u^\ft=\pi_0(h^\ft)\circ u^\ft=u.\]
We can then state the main result as follows:

\begin{theorem}\label{thm.main}
Let $A_K$ be an abelian variety over $K$. Then the  homomorphism
\[  {\rm H}^1(K,A_K)\to  {\rm Hom} (\pi_1(\GGr(A^\prime)), \Q/\Z)\]
mapping the torsor $X_K$ to the  homomorphism  
\[u_{X_K}\colon \pi_1(\GGr(A'))\to \Q/\Z \] in \eqref{eq.u} coincides with the homomorphism mapping $X_K$ to the homomorphism $v_{X_K}$ in \eqref{eq.v}.
If  furthermore the characteristic of $K$ is zero, then both constructions coincide with 
B\'egueri's construction in \eqref{eq.utau}, {\emph{i.e.}}, they explicate Shafarevich duality.
\end{theorem}

We expect that the construction  of \eqref{eq.v} via relative Picard functors can be related with that in \cite{To}.
 
\subsection{Equal characteristic case}

For $K$ of characteristic $p$, it is not clear either that the homomorphisms  in  \eqref{eq.u}, \eqref{eq.v} are still isomorphisms or that they provide Shafarevich's duality (cf.  \cite{Bes} for the good reduction case and in \cite{Ber} in general). 
However, we have a partial result on the prime-to-$p$ parts where Shafarevich's duality is quite easy to describe.

\subsubsection{Shafarevich's duality on the prime-to-$p$ parts} 
Let $n=l^r$ be a positive integer, prime to $p$, and large enough to kill the $l$-primary parts of the component groups of $A_K$ and $A_K'$.
Consider the perfect  cup product pairing \[\langle\ ,\ \rangle\colon {\rm H}^1(K,{}_nA_K)\times {}_nA^\prime_K(K)\to {\rm H}^1(K,\mu_n)=\Z/n\Z\]
on the(\'etale or flat) cohomology groups of the $n$-torsion points of $A_K$ and $A_K^\prime$.
Given an extension $\eta_n$ as in \eqref{eq.etan} (which corresponds to the torsor $X_K$) and  a point $a\in {}_nA^\prime_K(K)$; 
then $\langle \eta_n,a \rangle$ is the class of the pull-back along $a\colon \Z\to  {}_nA^\prime_K$ of the Cartier dual of $\eta_n$,
\[  \eta_n^D\colon 0\to \mu_n\to{}_nE_K^D \to {}_n A_K^\prime\to 0, \]
and it corresponds to the image of $a$ along the boundary map $ \partial\colon {}_n A_K^\prime(K) \to {\rm H}^1(K,\mu_n)$.
Furthermore, if  $  {}_nA^{\prime 0}$ denotes the quasi-finite subgroup of $n$-torsion sections of $A^{\prime 0}$,  we have 
\[\pi_1(\GGr(A^\prime))/n\cdot\pi_1(\GGr(A^\prime))=   {}_nA^{\prime 0}(R)={}_{n^2}A^{\prime}(R)/ {}_nA^{\prime}(R), \]
\[{}_n{\rm H}^1(K,A_K) = {\rm H}^1(K,{}_{n^2}A_K)/{\rm H}^1(K,{}_nA_K)  \]
(cf. \cite{Ber} \S1 )
and Shafarevich duality on the $d$-primary parts 
\begin{equation}\label{eq.shad}{}_n{\rm H}^1(K,A_K) \times
\pi_1(\GGr(A^\prime))/n\cdot\pi_1(\GGr(A^\prime))  \to {\rm H}^1(K,\mu_n)=\Z/n\Z \end{equation}
is induced by the above cup product.

The map $u^\tau\colon \pi_1(\GGr(A^\prime))\to \pi_0(\uuline{\rm H}^1(K,\mu_n))={\rm H}^1(K,\mu_n)$ associated to the torsor $X_K$ can also be viewed as the composition
\begin{equation}\label{eq.pi1cup}
\pi_1(\GGr(A^\prime))\stackrel{\delta}{\longrightarrow} {}_nA^\prime(R)={}_nA^\prime_K(K)\stackrel{\partial}{\longrightarrow} {\rm H}^1(K,\mu_n)=\Z/n\Z \subset \Q/\Z\end{equation}
where the first map is deduced from 
the exact sequence \begin{equation*} 
0\to {}_nA^\prime_K \to A^\prime_K\stackrel{n}{\to}  A^\prime_K \to 0 \end{equation*}
on passing to N\'eron models. More precisely we have 
\begin{equation*}
0\to {}_nA^\prime \to A^\prime\stackrel{n}{\to} nA^\prime \to 0 \end{equation*}
where $nA^\prime$ is a subgroup scheme of $A^\prime$ that contains $\tilde A^{\prime 0}$.
In particular, on applying the perfection of the Greenberg realization functor we get a homomorphism
\begin{equation} \label{eq.pi1A}\pi_1(\GGr(A^\prime) )=\pi_1(\GGr( nA^\prime ))\to  \pi_0({}_nA^\prime )={}_nA^\prime (R).
\end{equation}

\subsubsection{Comparison result}
Let $X_K$ be a torsor under $A_K$ of order $d$ with $d$  a power of a prime integer $l$, $l\neq p$.  
Let $n=l^r$ be a multiple of $d$ large enough to kill the  $l$-primary parts of the component groups of $A_K$ and $A_K'$. 
Fix an extension as \eqref{eq.eta} corresponding to $X_K$ and let $x_K=\spec(K')$ be a rigidificator of $\Pic_{X_K/K}$ contained in ${}_nE_K$.  
We show  that the composition of the maps in \eqref{eq.pi1cup} coincides with the map $u$ in \eqref{eq.u}. This is  sufficient to conclude that our construction via rigidificators (or equivalently via the relative Picard functor) is Shafarevich's duality on the prime-to-$p$ parts.

With notations as in \eqref{eq.eta}, observe that the $n$-multiplication on $A_K$ factors through $E_K$ so that we have a homomorphism $\gamma \colon E_K\to A_K$, with kernel ${}_n E_K$ such that $\gamma\circ \alpha=n$.  
Consider the sequence in 
\eqref{eq.muvEx}.
 We  have a diagram with exact rows
\begin{equation*}
\xymatrix{
0\ar[r]& {\mu_n} \ar[r]\ar@{=}[d] &  {}_nE_K^D \ar[r]\ar[d]&  \underline{\rm Ext}^1(A_K,\G_m) \ar[r]^{~ ~ n} \ar[d]^{\gamma^*}&
 A^\prime_K\ar[r]\ar@{=}[d]& 0\\
0\ar[r]& {\mu_n}=\underline{\rm Hom}(E_K,\G_m) \ar[r]& V_{x_K}^*\ar[r]& \underline{\rm Ext}^1(E_K,\G_m)_{x_K}\ar[r]&
 A^\prime_K\ar[r]& 0 .
}
\end{equation*}
Indeed   ${}_nE_K^D=\underline{\rm Hom}({}_nE_K,\G_m) $ maps canonically to $ V_{x_K}^*=\underline{\rm Mor}(x_K,\G_{m,K})$;  hence  ${}_nA^\prime_K$ maps to the torus $T_K=V_{x_K}^*/\mu_n$ in \eqref{eq.TExtA}. The push-out of the exact sequence $0\to {}_nA^\prime_K\to A^\prime_K\to A^\prime_K\to 0$  along   ${}_nA^\prime_K\to T_K$  provides the sequence \eqref{eq.TExtA}
and the homomorphism $\gamma^*$ sends a $\G_m$-extension of $A_K$ to its pull-back along $\gamma$ endowed with its canonical trivialization along $x_K$, induced by the canonical trivialization along   ${}_nE_K$. 

Moreover, the boundary map $ \partial\colon {}_n A_K'(K) \to {\rm H}^1(K,\mu_n)$ (of finite groups) is the composition of $ \nu\colon {}_n A_K^\prime(K)\to T_K(K)$ with the boundary map $h\colon T_K(K)\to {\rm H}^1(K,\mu_n)$, i.e. 
\[  \partial=h\circ \nu. \]
Recall furthermore that the kernel of the $n$-multiplication on $A'$ is a quasi-finite group scheme over $R$ whose finite part is an \'etale finite group scheme over $R$ of order prime to $p$, hence constant, because $R$ is strictly Henselian.

On the level of proalgebraic groups we then have a diagram with exact rows
\begin{equation*}
\xymatrix{
0\ar[r]& {}_n A^\prime(R)\ar[r]\ar[d]^{ \nu}& \GGr(A^\prime) \ar[r]\ar[d]^{\alpha^*}&
 \GGr(nA^\prime)  \ar[d]\ar[r]& 0\\
0\ar[r] & \GGr(T) \ar[r]& \GGr(j_* \underline{\rm Ext}^1(E_K,\G_m)_{x_K}) \ar[r]&
\GGr( A^\prime) \ar[r]& 0 
}
\end{equation*}
Since the vertical map on the left factors through a map $\nu^\ft\colon   {}_n A^\prime(R)\to  \GGr(T^\ft)$,  the homomorphism  
$u^\ft\colon \pi_1(\GGr (A^\prime))\to \pi_0 ( \GGr(T^\ft)  )=\pi_0(\GGr(T))_\tor $ in \eqref{eq.u} factors through the map $\delta\colon \pi_1(\GGr(A'))\to {}_n A^\prime(R)$ in \eqref{eq.pi1A} and hence
\[u^\tau= \partial\circ \delta=\pi_0(h^\ft)\circ \pi_0(\nu^\ft)\circ \delta=\pi_0(h^\ft)\circ u^\ft=u, \]
i.e.,  the homomorphism
$u\colon  \pi_1(\GGr (A^\prime))\to \Z/n\Z$ in \eqref{eq.u} coincides with that in  \eqref{eq.pi1cup}.
We conclude then
\begin{proposition}\label{pro.shap}
For any local field $K$ with algebraically closed residue field,  Shafarevich pairing coincides with the pairing constructed in Section \ref{sec.rig}, on the prime-to-$p$ parts.
\end{proposition}

The comparison for the $p$ parts is still open.
% and we think it should use  the PhD thesis \cite{Loe}.

  \thanks{\emph{Acknowledgements: } 
We thank M. Raynaud for some precious suggestions and J. Tong for pointing out some mistakes in the first draft. We thank Progetto di Eccellenza Cariparo 2008-2009 ``Differential methods in Arithmetics, Geometry and Algebra'' for a financial support.}

\end{document}